\documentclass[a4paper,12pt]{article}

\usepackage{amsmath, amssymb, amsthm}
\usepackage[dvips,final]{epsfig}
\usepackage[dvips,final]{graphicx}

\usepackage{color}

\definecolor{vio}{rgb}{0.5,0,0.5}
\definecolor{gre}{rgb}{0.1,0.6,0}
\definecolor{ora}{rgb}{0.8,0.2,0.1}

\usepackage{titlesec}
\titleformat{\section}{\bfseries}{\thesection}{1em}{}
\titleformat{\subsection}{\itshape}{\thesubsection}{1em}{}

\numberwithin{equation}{section}

\def\dd{\,\mathrm{d}}
\def\dive{\mathrm{\,div\,}}
\def\pif{\hbox{\ if}\ }
\def\ejt{e_j(t)}
\def\sumjm{\sum_{j=-m}^m}
\def\sumkm{\sum_{k=1}^m}
\def\sumlm{\sum_{l=0}^m}

\def\real{\mathbb{R}}
\def\nat{\mathbb{N}}

\def\AA{\mathbf{A}}
\def\BB{\mathbf{B}}
\def\play{\mathfrak{p}}

\def\nas{\nabla_s}

\def\io{\int_{\Omega}}
\def\ippi{\int_{2\pi}^{4\pi}}
\def\ipo{\int_{\partial\Omega}}

\def\om{^{(m)}}
\def\os{^{(\tau)}}
\def\pp{_{2\pi}}

\def\supess{\mathop{\mathrm{sup\,ess}}}

\def\be{\begin{equation}\label}
\def\ee{\end{equation}}

\def\ber{\begin{eqnarray}}
\def\eer{\end{eqnarray}}

\def\bers{\begin{eqnarray*}}
\def\eers{\end{eqnarray*}}

\newfont{\ctv}{msam10}

\newcommand{\bbox}{\mbox{\ctv \symbol{4}}}
\def\QED{{${}\hfill\bbox$}}
\newenvironment{pf}[1]{\par\vskip1mm{\noindent\it #1.}\ }{\QED\par
\vskip2mm}

\def\bpf{\begin{pf}}
\def\epf{\end{pf}}

\newtheorem{theorem}{Theorem}[section]

\newtheorem{hypothesis}[theorem]{Hypothesis}
\newtheorem{proposition}[theorem]{Proposition}

\begin{document}

\title{Periodic waves in unsaturated porous media with hysteresis 
\thanks{Supported by the GA\v CR Grant GA15-12227S, RVO: 67985840, FWF P23628-N18,
and by the FP7-IDEAS-ERC-StG \#256872 (EntroPhase)}.}

\author{Bettina Albers
\thanks{University of Duisburg-Essen, 
Faculty of Engineering, Department of Civil Engineering,
45117 Essen , Germany, E-mail {\tt bettina.albers@uni-due.de}.}
\and Pavel Krej\v c\'{\i}
\thanks{Institute of Mathematics, Czech Academy of Sciences, \v Zitn\'a~25,
CZ-11567~Praha 1, Czech Republic, E-mail {\tt krejci@math.cas.cz}.}
\and Elisabetta Rocca
\thanks{Dipartimento di Matematica, Universit\`a degli Studi di Pavia. Via Ferrata 5, 
I-27100 Pavia, Italy, E-mail {\tt elisabetta.rocca@unipv.it}.}
}

        
\maketitle


\begin{abstract}
We consider a PDE system with degenerate hysteresis describing unsaturated flow
in 3D porous media. Assuming that a time periodic forcing is prescribed
on the boundary, we prove that a time periodic response exists as long as
the amplitude of the forcing terms is small enough to keep the solution within
the convexity domain of the hysteresis operator.
\end{abstract}


\section*{Introduction}\label{int}

Periodic waves of a given frequency are often used in non-destructive testing of
porous media. In particular, building materials, geomaterials, tissues, nanomaterials etc. are examples
of porous materials, in which non-destructive testing methods are of central importance.
The reason is not only that the sample does not have to be destroyed or invaded, but another
advantage in comparison with conventional methods is that non-destructive testing is mostly cost-saving.
For example, this is the case of testing of soils, where procedures exploiting the
properties of acoustic waves are cheaper than drilling the boreholes.

There are many different techniques of non-destructive testing: 
ultrasonic methods, magnetic particle inspection, liquid penetrant inspection,
electrical measurements or radiography -- for details see, e.~g., \cite{Shull, Veenstra, CosentiniFoti}.
For further information on acoustic methods of non-destructive testing see, e.~g., \cite{Gan,Kraut}.

In the geotechnical field, the wave analysis of both body and surface waves (\cite{Hab})
may lead to the construction of several non-destructive testing methods. Body or bulk waves travel
through the interior of a medium, while surface waves propagate along the surface of a body or along
the interface of two media. The amplitudes of surface waves decay in the direction perpendicular
to the surface so fast that they can be assumed to be zero in the depth of a few wavelengths.
The analysis of surface waves in saturated porous media (\cite{a2006}) may help to develop
a method for soil characterization. By use of the SASW-technique (Spectral Analysis of Surface Waves),
e.~g., \cite{Roesset,Foti}, conclusions about building grounds can be drawn from the measurement
of sound wave speeds. In other words, expensive and invasive acoustic measurements in boreholes
or laboratory tests are not necessary to characterize the soil prior to a building project.
In SASW tests, two or more receivers are placed on the surface, and a hammer
(or a signal with a certain frequency) is used to generate surface waves whose speeds are recorded.
Algorithms based on the Fast Fourier Transform applied to the acquired data then produce
a stiffness versus depth plot. While in the classical method the wave propagation in single-component
media is analyzed, in \cite{a2006} the considerations are extended to two and three
component modeling.

If the pores of a porous medium are filled with two (or more) immiscible fluids as, for example,
water and air, then they are called ‘partially saturated’. The pore fluids possess different
partial pressures, i.~e. there exists a discontinuity in the pressure across the interface.
This difference is called the capillary pressure. It depends on the geometry of the pore space,
on the nature of the solids and on the degree of saturation, i.~e. the ratio of the volume
occupied by one of the pore fluids over the entire pore volume.

For the description of the propagation of sound waves in partially saturated soils (three-component media)
a linear macroscopic model is introduced in \cite{Hab}. 
However, experimental studies of wetting and dewetting curves in partially saturated porous media
exhibit strong capillary hysteresis effects, see \cite{hrrp}, which are due to the surface tension
on the liquid-gas interface. Flynn et al. \cite{Flynn, flpokr}
suggested to model hysteresis phenomena in porous media by means of the Preisach operator
originally designed in \cite{pr} for magnetic hysteresis in ferromagnetics.
Two models to describe processes
in partially saturated media are presented in \cite{ak2016}. The first model does not
explicitly contain a hysteresis operator and the effect of hysteresis in the capillary pressure
curve is accounted for by investigating the two processes drainage and imbibition separately,
cf.~also \cite{CEER,Murphys}. The second model is a thermomechanical model involving the Preisach operator,
while plastic hysteresis is described in terms of the Prandtl-Reuss model.

In the present paper, we study the mathematical problem of well posedness of the porous medium
model proposed in \cite{ak} under periodic mechanical forcing.
We assume that a time periodic force is prescribed on the boundary of the domain, and
look for time periodic mechanical waves in the system of balance equations.
The main difference with respect to \cite{ak} consists in the hypothesis
that the solid matrix material is elastic within the small deformation hypothesis,
so that the momentum balance equation is linear. The only nonlinearity in the problem
is thus the degenerate Preisach hysteresis operator in the mass balance equation. On the other hand,
since viscosity is missing in the model, we lose the higher order a priori estimates,
which were used in \cite{ak} to control the degeneracy of the Preisach operator.
Instead, we make use of the second order energy inequality related to the convexity
of small amplitude hysteresis loops to prove that periodic solutions of the system
exist provided the amplitude of the external forcing is sufficiently small.
Note that for any nonlinear pressure-saturation relation without hysteresis,
such a result would be much more difficult to obtain, since no counterpart
of the second order energy inequality is available in this case.

The structure of the paper is as follows. In Section \ref{mod}, we present the model
situation, and in Section~\ref{sta} we state Theorem~\ref{t1} which is the main Existence
Theorem of the paper. Section~\ref{hys} is devoted to a survey about Preisach hysteresis,
and Section \ref{reg} contains the proof of Theorem~\ref{t1}.


\section{The model}\label{mod}

\begin{figure}[htb]
\centerline{\psfig{file=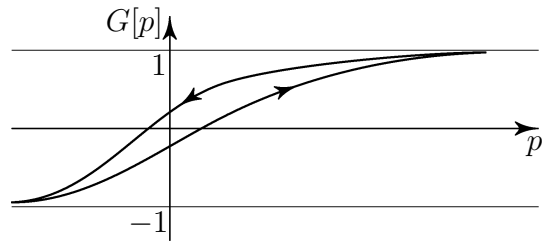,width=7.5cm}}
\caption{The pressure-saturation hysteresis}
\label{f2}
\end{figure}

The present paper deals with the following model for fluid flow in an unsaturated porous solid
\ber\label{e1}
\rho_S u_{tt} + c u_t &=& \dive \AA\nas u + \nabla p + f_0\,,\\ \label{e2}
G[p]_t &=& \dive u_t + \frac{\mu}{\rho_L} \Delta p\,.
\eer
in a spatial domain $\Omega$ and time $t \in \real$, for unknown functions
$u$ (displacement) and $p$ (capillary pressure), with a operator $G$
characterizing the hysteresis dependence between $p$ and the relative air content $A \in [-1,1]$, 
$A = G[p]$ as in \cite{ak}, see Figure \ref{f2}, and with constant coefficients
$\rho_S, \rho_L$ (mass densities of the solid and liquid, respectively), and $\mu$ (permeability).
By $\nas$ we denote the symmetric gradient, $\AA$ is a constant elasticity matrix,
and $f_0$ is a given external volume force. Eq.~\eqref{e1} is the momentum balance,
\eqref{e2} is the liquid mass balance.

On the boundary $\partial\Omega$ we prescribe boundary conditions
\be{bc}
(u - u^*)\big|_{\partial\Omega} = 0\,, \quad \nabla p \cdot n\big|_{\partial\Omega} = \gamma(x)(p^* - p)\,,
\ee
where $u^*$ is a given displacement, $n$ is the unit outward normal vector,
$p^* = p^*(x,t)$ is a given outer pressure, and $\gamma(x) \ge 0$ is a given permeability
of the boundary which is positive on a set of positive measure.

A similar system was derived in \cite{ak} in the form
\ber\label{e1t}
\rho_S u_{tt} &=& \dive (\BB\nas u_t + P[\nas u]) + \nabla p + f_0\,,\\ \label{e2t}
G[p]_t &=& \dive u_t + \frac{1}{\rho_L} \dive(\mu(p)\nabla p)\,.
\eer
as a model for isothermal flow in an unsaturated viscoelastoplastic porous solid, where
$P$ is a constitutive operator of elastoplasticity, $\BB$ is a constant viscosity matrix, 
and $\mu(p)$ is a pressure dependent permeability coefficient. In reality, the coefficient
$\mu$ should depend on the saturation, that is, on $G[p]$, but the analysis carried out
in \cite{bv1,bv2} shows that the presence of the hysteresis operator in the coefficient
makes the problem difficult, so that it cannot be solved without additional space or time regularization.

The main result of \cite{ak} was the proof of existence of a strong solution to the above
system with natural boundary conditions and given initial conditions.
In particular, it was shown that the solution remains bounded away from the degeneracy
of the hysteresis operator $G$.

The system \eqref{e1}--\eqref{e2} differs from \eqref{e1t}--\eqref{e2t} in several respects.
Notice first that the strong viscous dissipative term $\dive \BB \nas u_t$ is replaced
with a more realistic weaker term $c u_t$ corresponding to contact friction
on the solid-liquid interface with a constant friction parameter $c>0$.
The strong a priori estimates as in \cite{ak} resulting from the viscous term
are no longer available for \eqref{e1}--\eqref{e2}, so that we need additional modeling hypotheses
to prove the solvability of the system. More specifically,
the solid matrix is assumed elastic, that is, $P[\nas u] = \AA \nas u$, 
the permeability coefficient $\mu$ is a positive constant, and we consider
time periodic data.

The system is linear in $u$, so that we can replace $u$ by $u - u^*$ if $u^*$ is sufficiently
regular, and reformulate \eqref{e1}--\eqref{bc} as
\ber\label{e1a}
\rho_S u_{tt} + c u_t &=& \dive \AA\nas u + \nabla p + f\,,\\ \label{e2a}
G[p]_t &=& \dive u_t + \frac{\mu}{\rho_L} \Delta p + h\,,
\eer
\be{bca}
u \big|_{\partial\Omega} = 0\,, \quad \nabla p \cdot n\big|_{\partial\Omega} = \gamma(x)(p^* - p)\,,
\ee
with given functions $f$ and $h$.


\section{Statement of the problem}\label{sta}

We assume that all given functions $f$, $h$, and $p^*$ are time periodic with the same period.
The values of the physical constants are not relevant for our analysis, so that for simplicity,
we consider system \eqref{e1a}--\eqref{bca} in the form
\ber\label{e1b}
u_{tt} + u_t &=& \dive \AA\nas u + \nabla p + f\,,\\ \label{e2b}
G[p]_t &=& \dive u_t + \Delta p + h\,,
\eer
\be{bcb}
u \big|_{\partial\Omega} = 0\,, \quad \nabla p \cdot n\big|_{\partial\Omega} = \gamma(x)(p^* - p)
\ee
with $2\pi$-periodic data, and introduce the notation
\be{lp}
L^q\pp (\Omega) = \{y \in L^q_{loc}(\Omega\times \real): y(x,t+2\pi) = y(x,t)\ \mbox{a.\,e.}\},
\ee
and similarly for $L^q\pp(\partial\Omega)$ etc. The norm in $L^q\pp (\Omega)$ is defined as
\be{nlp}
\|y\|_{q,\Omega,2\pi} = \left(\int_{2\pi}^{4\pi}\io |y|^q \dd x\dd t\right)^{1/q},
\ee
and in $L^q\pp(\partial\Omega)$ we introduce the seminorm
\be{nlpb}
\|y\|_{q,\partial\Omega,2\pi,\gamma}
= \left(\int_{2\pi}^{4\pi}\ipo \gamma(x)|y|^q \dd s(x)\dd t\right)^{1/q}.
\ee

\begin{theorem}\label{t1}
Let $\Omega$ be a bounded domain with $C^{1,1}$ boundary, let $\gamma \in C^1(\partial\Omega)$
be a nonnegative function which does not identically vanish, and let the data $f, h, p^*$ be such that
$f, f_t, h, h_t \in L^2\pp(\Omega)$, $p^*, p^*_t \in L^2\pp(\partial\Omega)$.
Set
$$
\delta = \max\{\|f\|_{2,\Omega,2\pi}, \|f_t\|_{2,\Omega,2\pi}, \|h\|_{2,\Omega,2\pi}, \|h_t\|_{2,\Omega,2\pi},
\|p^*\|_{2,\partial\Omega,2\pi,\gamma}, \|p^*_t\|_{2,\partial\Omega,2\pi,\gamma}\}.
$$
Then there exists $\delta^* > 0$ such that if $\delta < \delta^*$, then
system \eqref{e1b}--\eqref{bcb} has a solution $u,p$ such that
$u, u_t, u_{tt}, \nas u, \nas u_t, \dive \AA \nas u, p, p_t, \nabla p, \nabla p_t, \Delta p
\in L^2\pp(\Omega)$.
\end{theorem}

The reason why we have to assume that the data are small is related to the fact that
higher order a priori estimates, which are not available here due to the absence
of the viscosity term $\dive\BB\nas u_t$, can only be recovered as long as the input
$p$ of the hysteresis operator $G$ stays in the convexity domain of $G$. Details
will be given in Section \ref{hys}.


\section{Hysteresis operators}\label{hys}

We recall here the basic concepts of the theory of hysteresis operators that are needed in the sequel.
The construction of the operator $G$ is based on the variational inequality
\be{play}
\left\{
\begin{array}{ll}
|p(t) - \xi_r(t)| \le r & \forall t\in [0,T]\,,\\
(\xi_r(t))_t(p(t) - \xi_r(t) - z) \ge 0 & \mbox{a.~e.}\ \forall z \in [-r,r]\,,\\
p(0) - \xi_r(0) = \max\{-r, \min\{p(0), r\}\}.
\end{array}
\right.
\ee
It is well known (\cite{kp}) that for each given input function $p \in W^{1,1}(0,T)$
for some $T>0$ and each
parameter $r>0$, there exists a unique solution $\xi_r\in W^{1,1}(0,T)$ of the variational inequality
\eqref{play}.
The mapping $\play_r: W^{1,1}(0,T) \to W^{1,1}(0,T)$ which with each $p\in W^{1,1}(0,T)$ associates
the solution $\xi_r = \play_r[p] \in W^{1,1}(0,T)$ of \eqref{play} is called
the {\em play operator\/}, and the parameter
$r>0$ can be interpreted as a {\em memory parameter\/}.
The proof of the following statements
can be found, e.~g., in \cite[Chapter II]{book}.

\begin{proposition}\label{p5}
For each $r>0$, the mapping $\play_r: W^{1,1}(0,T) \to W^{1,1}(0,T)$ is Lipschitz continuous
and admits a Lipschitz continuous extension to $\play_r: C[0,T] \to C[0,T]$ in the sense that
for every $p_1, p_2 \in C[0,T]$ and every $t \in [0,T]$ we have
\be{lipc}
|\play_r[p_1](t) - \play_r[p_2](t)| \le \|p_1 - p_2\|_{[0,t]} := \max_{\tau \in [0,t]}|p_1(\tau) - p_2(\tau)|\,.
\ee
Moreover, for each $p \in W^{1,1}(0,T)$, the energy balance equation
\be{enerpl}
\play_r[p]_t p - \frac12\left(\play_r^2[p]\right)_t = \left|r\, \play_r[p]_t\right|
\ee
and the identity
\be{mono}
\play_r[p]_t p_t = (\play_r[p]_t)^2
\ee
hold almost everywhere in $(0,T)$.
\end{proposition}

Similarly as above, we define the spaces of $2\pi$-periodic functions of time
$$
L^{q}\pp = \{z \in L^{q}_{loc}(\real): z(t+2\pi) = z(t)\ a.~e.\}, \ 
C\pp =  \{z \in C(\real): z(t+2\pi) = z(t)\ \forall t\in \real\}\,,
$$
endowed with the natural norms
$$
|z|_{q,2\pi} = \left(\int_{2\pi}^{4\pi} |z(t)|^q\dd t\right)^{1/q}
$$
and similarly for $C\pp, W^{k,q}\pp$ for $k \in \nat$ and $q\ge 1$, etc.

\begin{proposition}\label{p6}
For every $p \in C\pp$ and every $r>0$ we have $\play_r[p](t+2\pi) = \play_r[p](t)$
for all $t \ge 2\pi$. In particular, by extending $\play_r[p]$ backward from the interval
$[2\pi, \infty)$ periodically to $\real$, we can assume that $\play_r$ maps $C\pp$ into $C\pp$.
\end{proposition}

Given a nonnegative function $\rho \in L^1((0,\infty)\times \real)$,
we define the operator $G$ as a mapping that with each $p \in C\pp$ associates the integral
\be{pre}
G[p](t) = \int_0^\infty \int_0^{\play_r[p](t)} \rho (r,v)\dd v\dd r\,.
\ee
Directly from the definition \eqref{play} of the play, we see with the notation of \eqref{lipc} that
the implication
\be{null}
r\ge \|p\|_{[0,t]} \ \Longrightarrow \ \play_r[p](t) = 0
\ee
holds for every $p \in W^{1,1}(0,T)$ and every $T>0$ (hence, for every $p \in C\pp$),
so that the integration domain in \eqref{pre} is always bounded.

Definition \eqref{pre} is equivalent to the {\em Preisach model\/} proposed in \cite{pr}, see \cite{max}. 
For our purposes, we prescribe the following hypotheses on $\rho$.

\begin{figure}[htb]
\centerline{\psfig{file=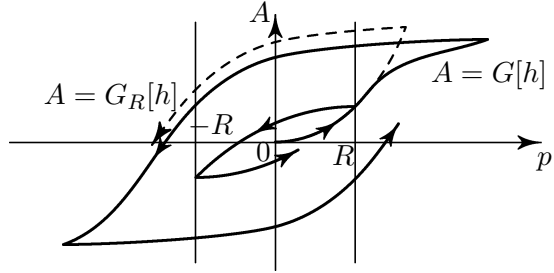,width=7.5cm}}
\caption{Local behavior of the operator $G$ and its convexification $G_R$.}
\label{f1}
\end{figure}

\begin{hypothesis}\label{h2}
The function $\rho \in W^{1,\infty}((0,\infty)\times \real)$ is such that
there exists a function $\rho^* \in L^1(0,\infty)$ such that for a.~e.
$(r,v) \in (0,\infty)\times\real$ we have $0 \le \rho(r,v) \le \rho^*(r)$, and we put
\be{crho}
C_\rho = \int_0^\infty\int_{-\infty}^\infty\rho(r,v)\dd v\dd r\,, \quad
C_\rho^* = \int_0^\infty\rho^*(r)\dd r\,.
\ee
Furthermore, there exists $R>0$ for which the following condition holds:
\be{con}
A_R := \inf\{\rho(r,v): r+|v| \le R\} > 0\,.
\ee
\end{hypothesis}
Put $C_R := \sup\{\left|\frac{\partial}{\partial v} \rho(r,v)\right|: r+|v| \le R\}$. Taking
$R>0$ smaller, if necessary, we can assume that there exists $K_R>0$ such that
\be{con2}
\frac12 A_R - R C_R \ge K_R\,.
\ee
{}From \eqref{enerpl}, \eqref{mono}, and \eqref{pre} we immediately deduce the Preisach
energy identity
\be{enerpr}
G[p]_t p - V[p]_t = |D[p]_t| \ \mbox{ a.~e.}
\ee
with a Preisach potential $V$ and dissipation operator $D$ defined as
\be{prpot}
V[p](t) = \int_0^\infty \int_0^{\play_r[p](t)} v \rho (r,v)\dd v\dd r\,, \quad
D[p](t) = \int_0^\infty \int_0^{\play_r[p](t)} r \rho (r,v)\dd v\dd r\,.
\ee
A straightforward computation shows that $G$ is Lipschitz continuous
in $C[0,T]$. Indeed, using \eqref{lipc} and Hypothesis \ref{h2}, we obtain for $p_1, p_2 \in C[0,T]$ and
$t \in [0,T]$ that
\be{lipg}
|G[p_1](t) - G[p_1](t)| = \left|\int_0^\infty\int_{\play_r[p_1](t)}^{\play_r[p_2](t)}\rho(v,r)\dd v\dd r\right|
\le C_\rho^* \max_{\tau \in [0,t]} |p_1(\tau) - p_2(\tau)|\,.
\ee
Following \cite{mhd}, we define the convexified operator $G_R$ (see Figure \ref{f1}) by a formula similar to \eqref{pre}
\be{prec}
G_R[p](t) = \int_0^\infty \int_0^{\play_r[p](t)} \rho_R (r,v)\dd v\dd r\,,
\ee
where
\begin{equation}\label{psiR}
\rho_R(r,v) = \left \{
\begin{array}{ll}
 \rho(r,v) \qquad & \pif r+|v| \le R,\\[2mm]
 \rho(r, -R+r) \qquad & \pif v < - R + r, \, r \le \, R,\\[2mm]
 \rho(r, R - r) \qquad & \pif v > R - r, \, r \le R,\\[2mm]
 \rho(R, 0), \qquad & \pif r > R.
\end{array}
\right.
\end{equation}

It is shown in \cite{mhd} that the operator $G_R$ satisfies globally
the hypotheses of \cite[Theorem II.4.19]{book}, that is, the ascending hysteresis branches are
uniformly convex and the descending branches are uniformly concave, so that for each input
$p \in W^{2,1}\pp$, the second order energy inequality
holds in the form
\be{ene2}
-\ippi G_R[p]_{t} p_{tt} \dd t \ge \frac{K_R}{2} \ippi |p_t|^3 \dd t\,.
\ee
As a consequence of \eqref{enerpr} we also have
\be{ene3}
\ippi G_R[p]_{t} p \dd t \ge 0
\ee
for every $p \in W^{1,1}\pp$.
Since the density $\rho$ is globally bounded above by a constant $H_\rho$ and the implication
\eqref{null} holds, the operator $G_R$ has quadratic growth in the sense
\be{qua}
|G_R[p](t)| \le H_\rho \|p\|_{[0,t]}^2\,, \quad |G_R[p]_t(t)| \le H_\rho \|p\|_{[0,t]} |p_t(t)| \ \mbox{a.~e.}
\ee
for all functions $p \in W^{1,1}(0,T)$ and all $T>0$.
On the other hand, we have the implication
\be{conve}
|p(t)|\le R \ \forall t\ge 0 \ \Longrightarrow \ G_R[p](t) = G[p](t) \ \forall t\ge 0\,.
\ee


\section{Proof of Theorem \ref{t1}}\label{reg}

We replace the operator $G$ by $G_R$ and consider the variational formulation
of the convexified version of Problem \eqref{e1b}--\eqref{bcb}
\ber\label{e1r}
\io ((u_{tt} + u_t)\phi + \AA\nas u :\nas\phi + p \dive\phi)\dd x
&=& \io f\phi \dd x\,,\\ \label{e2r}
\io ((G_R[p]_t - \dive u_t)\psi + \nabla p \nabla\psi)\dd x
&=& \io h\psi \dd x + \ipo \gamma(x)(p^* - p)\psi\dd s(x)
\eer
for every test functions $\phi \in W^{1,2}_0(\Omega;\real^3)$ and $\psi \in W^{1,2}(\Omega)$.

With the intention to use the Galerkin method, we choose $\{\phi_k; k=1,2,\dots\}$
in $L^2(\Omega;\real^3)$ and $\{\psi_l; l=0,1,2,\dots\}$ in $L^2(\Omega)$
to be the complete orthonormal systems of eigenfunctions defined by
\be{eigen}
-\dive\AA \nas \phi_k = \lambda_k \phi_k \ \text{ in } \ \Omega\,, \ \ \phi_k\big|_{\partial\Omega} = 0\,,
\quad -\Delta \psi_l = \mu_l \psi_l \ \text{ in } \ \Omega\,,
\ \ \nabla \psi_l\cdot n\big|_{\partial\Omega} = 0\,,
\ee
with $\mu_0 = 0$, $\lambda_k>0, \mu_l > 0$ for $k,l \ge 1$.


\subsection{Galerkin approximations}\label{gale}

Approximate $2\pi$-periodic solutions
will be searched in the form
\be{app}
u\om (x,t) = \sumjm\sumkm u_{jk}\ejt\phi_k(x)\,, \quad p\om(x,t) = \sumjm\sumlm p_{jl}\ejt\psi_l(x)\,,
\ee
with
$$
\ejt = 
\left\{
\begin{array}{ll}
\sin jt & \pif\ j \ge 1\,,\\
\cos jt & \pif\ j \le 0\,, 
\end{array}\right.
$$
and with real coefficients $u_{jk}, p_{jl}$ which satisfy the system
\ber \nonumber
&&\ippi\io ((u\om_{tt} + u\om_t)\phi_k(x) + \AA\nas u\om :\nas\phi_k(x) + p\om \dive\phi_k(x))\ejt\dd x \dd t\\  \label{e1g}
&& \qquad = \ippi\io f\phi_k(x)\ejt \dd x\dd t\,,\\ \nonumber
&&\ippi\io ((G_R[p\om]_t - \dive u\om_t)\psi_l(x) + \nabla p\om \nabla\psi_l(x))\ejt\dd x \dd t\\ \label{e2g}
&&\qquad  =\ippi\io h\psi_l(x)\ejt \dd x \dd t + \ippi\ipo \gamma(x)(p^* - p\om)\psi_l(x)\ejt\dd s(x) \dd t\,.
\eer
for $j=-m, \dots, m$, $k=1, \dots, m$, $l=0,\dots, m$. This is an algebraic problem of
$(2m+1)^2$ equations for a $(2m+1)^2$-dimensional real unknown vector
\be{U}
U = (u_{jk}, p_{jl}),\ j=-m, \dots, m,\ k=1, \dots, m,\ l=0,\dots, m
\ee
that we solve by using the degree theory. We define a continuous family of mappings
$\mathcal{T}_\alpha: \real^{(2m+1)^2} \to \real^{(2m+1)^2}$ for $\alpha \in [0,1]$ which
with $U$ as in \eqref{U} associate the vector $\mathcal{T}_\alpha(U) = V$ of the form
\be{V}
V = (v_{jk}, w_{jl}),\ j=-m, \dots, m,\ k=1, \dots, m,\ l=0,\dots, m,
\ee
given by the formula
\ber \nonumber
v_{jk}&=&\ippi\io ((u\om_{tt} + u\om_t)\phi_k(x) + \AA\nas u\om :\nas\phi_k(x) + p\om \dive\phi_k(x))\ejt\dd x \dd t\\  \label{V1}
&& - \ippi\io \alpha f\phi_k(x)\ejt \dd x\dd t\,,\\ \nonumber
w_{jl}&=&\ippi\io (((1-\alpha)p\om + \alpha G_R[p\om])_t - \dive u\om_t)\psi_l(x) + \nabla p\om \nabla\psi_l(x))\ejt\dd x \dd t\\ \label{V2}
&&  -\ippi\io \alpha h\psi_l(x)\ejt \dd x \dd t
- \ippi\ipo \gamma(x)(\alpha p^* - p\om)\psi_l(x)\ejt\dd s(x) \dd t\,.
\eer
System \eqref{e1g}--\eqref{e2g} can be interpreted as $\mathcal{T}_1(U) = 0$. Clearly, all mappings
$\mathcal{T}_\alpha$ for $\alpha \in [0,1]$ are continuous. We now show that the equation
$\mathcal{T}_\alpha(U) = 0$ for any $\alpha \in [0,1]$ has no solution $U$ on the boundary
of a sufficiently large ball $\mathcal{B}_K \subset \real^{(2m+1)^2}$.
Indeed, assume that $V = \mathcal{T}_\alpha(U) = 0$ for some $U$ and $\alpha$. Note that
\be{deri}
u\om_t(x,t) = \sumkm\sumjm \hat u_{jk}\ejt\phi_k(x)\ \ \mbox{with }\ \hat u_{jk} = j u_{-jk}\,.
\ee
We then have
\bers
0 &=& \sumkm\sumjm v_{jk} \hat u_{jk} + \sumlm\sumjm w_{jl} p_{jl}\\
&=& \ippi\io (|u\om_t|^2 - \alpha f u\om_t + \alpha G_R[p\om]_t p\om + |\nabla p\om|^2
- \alpha h p\om)\dd x\dd t\\
&& - \ippi\ipo \gamma(x)(\alpha p^* - p\om)p\om \dd s(x) \dd t\,,
\eers
and using \eqref{ene3} we obtain
\bers
&&\ippi\io (|u\om_t|^2  + |\nabla p\om|^2)\dd x\dd t + \ippi\ipo \gamma(x)|p\om|^2 \dd s(x) \dd t\\
&&\qquad \le \alpha \ippi\io (f u\om_t + h p\om)\dd x\dd t
+ \alpha\ippi\ipo \gamma(x)p^*p\om \dd s(x) \dd t\,,
\eers
so that
$$
 \sumkm\sumjm |u_{jk}|^2 + \sumlm\sumjm |p_{jl}|^2 \le K
$$
independently of $\alpha$. We see that $\mathcal{T}_\alpha$ is a homotopy of continuous mappings
on $\real^{(2m+1)^2}$ and such that the equation $\mathcal{T}_\alpha(U) = 0$ has no solution
on the boundary of any ball of radius bigger than $K$. Since $\mathcal{T}_0$ is odd,
its topological degree with respect to the ball $\mathcal{B}_{K+1}$ and the point $0$ is nonzero,
and remains constant for all $\alpha \in [0,1]$. We conclude that the equation $\mathcal{T}_1(U) = 0$
has a solution, which, by definition satisfies \eqref{e1g}--\eqref{e2g}, as well as the estimate
\be{es1}
\|u\om_t\|_{2,\Omega,2\pi} + \|\nabla p\om\|_{2,\Omega,2\pi}
+ \|p\om\|_{2,\partial\Omega,2\pi,\gamma} \le C\delta
\ee
with $\delta$ from Theorem \ref{t1} and with a constant $C$ independent of $m$ and $\delta$.

By iterating the formula \eqref{deri} we obtain
\bers
u\om_{ttt}(x,t) &=& \sumkm\sumjm u^\sharp_{jk}\ejt\phi_k(x)\ \ \mbox{with }\ u^\sharp_{jk} = -j^3 u_{-jk}\,,\\
p\om_{tt}(x,t) &=& \sumlm\sumjm p^\sharp_{jl}\ejt\psi_l(x)\ \ \mbox{with }\ p^\sharp_{jl} = -j^2 p_{jl}\,.
\eers
We now test \eqref{e1g} by $-u^\sharp_{jk}$, \eqref{e2g} by $-p^\sharp_{jl}$, and use \eqref{ene2} to obtain
\be{es2}
\|u\om_{tt}\|_{2,\Omega,2\pi}^2 + \|p\om_t\|_{3,\Omega,2\pi}^3
+ \|\nabla p\om_t\|_{2,\Omega,2\pi}^2 + \|p\om_t\|_{2,\partial\Omega,2\pi,\gamma}^2 \le C\delta^2
\ee
with a constant $C$ independent of $m$ and $\delta$. Then, testing \eqref{e1g} by $\bar u_{jk}
:= -j^2 u_{jk}$, we obtain using \eqref{es1}, \eqref{es2} that
\be{es3}
\|\nas u\om_{t}\|_{2,\Omega,2\pi} \le C\delta\,.
\ee
By \eqref{qua} we have
\be{e4}
|G_R[p\om]_t(x,t)| \le H_\rho \|p\om(x,\cdot)\|_{[2\pi, 4\pi]} |p\om_t(x,t)|
\ee
for all $x \in\Omega$ and $t \in [2\pi, 4\pi]$. Using the inequality
\be{e3}
\|p\om(x,\cdot)\|_{[2\pi, 4\pi]}\le\frac{1}{2\pi}\ippi |p\om(x,\tau)|\dd\tau+\ippi |p\om_t(x,\tau)|\dd\tau\,,
\ee
we have, by \eqref{e3} and Minkowski's inequality,
$$
\left(\io \|p\om(x,\cdot)\|_{[2\pi, 4\pi]}^6 \dd x\right)^{1/6}
\le C\ippi \left(\io (|p\om(x,t)|^6 + |p\om_t(x,t)|^6)\dd x\right)^{1/6}\dd t\,,
$$
so that, by the Sobolev embedding and estimates \eqref{es1}--\eqref{es2},
\bers
&&\hspace{-8mm} \left(\io \|p\om(x,\cdot)\|_{[2\pi, 4\pi]}^6 \dd x\right)^{1/6}\\
&\le& C\ippi \left(\io (|p\om(x,t)|^2 + |p\om_t(x,t)|^2 + |\nabla p\om(x,t)|^2
+ |\nabla p\om_t(x,t)|^2)\dd x\right)^{1/2} \dd t \le C\delta\,.
\eers
{}From \eqref{e4} and H\"older's inequality it follows
\be{e5}
\io|G_R[p\om]_t(x,t)|^2\dd x \le C \delta^2 \left(\io |p\om_t(x,t)|^3\dd x\right)^{2/3}\,,
\ee
hence
\be{es4}
\ippi\left(\io|G_R[p\om]_t(x,t)|^2\dd x\right)^{3/2}\dd t \le C \delta^5\,.
\ee


\subsection{Passage to the limit}\label{pass}

The compactness argument will be based on an anisotropic embedding formula which is
a special case of the theory developed in \cite{bin}.
For a bounded domain $D\subset \real^N$, an open bounded interval $\omega \subset \real$, and real numbers
$q,r \ge 1$ we define for $v \in L^r(\omega; L^q(D))$ and $w \in L^q(D; L^r(\omega))$ their
anisotropic norms
\be{lrq}
\|v\|_{q,r} = \left(\int_\omega \left(\int_D |v(x,t)|^q\dd x\right)^{r/q}\dd t\right)^{1/r},
\quad \|w\|^*_{r,q} = \left(\int_D \left(\int_\omega |w(x,t)|^r\dd t\right)^{q/r}\dd x\right)^{1/q}.
\ee
We also introduce the anisotropic Sobolev spaces
\bers
W^{r_0, q_0; r_1, q_1}(\omega, D) &=& \left\{v \in L^1(D\times \omega):
\frac{\partial v}{\partial t} \in  L^{r_0}(\omega; L^{q_0}(D))\,,
\ \nabla v \in L^{r_1}(\omega; L^{q_1}(D))\right\},\\
W^{q_0, r_0; q_1, r_1}(D, \omega) &=& \left\{w \in L^1(D\times \omega):
\frac{\partial w}{\partial t} \in  L^{q_0}(D; L^{r_0}(\omega))\,,
\ \nabla w \in L^{q_1}(D; L^{r_1}(\omega))\right\}.
\eers
We will repeatedly use the following compact embedding result.

\begin{proposition}\label{emb}
Let the domain $D$ be Lipschitzian.
\begin{itemize}
\item [{\rm (i)}] If $q \ge \max\{q_0, q_1\}$,  $r \ge \max\{r_0, r_1\}$, and
$$
\left(1 - \frac{1}{r_0} + \frac{1}{r}\right)\left(\frac{1}{N} - \frac{1}{q_1} + \frac{1}{q}\right)
> \left(\frac{1}{r_1} - \frac{1}{r}\right)\left(\frac{1}{q_0} - \frac{1}{q}\right),
$$
then $W^{r_0, q_0; r_1, q_1}(\omega, D)$ is compactly embedded in $L^r(\omega; L^q(D))$ and
 $W^{q_0, r_0; q_1, r_1}(D, \omega)$ is compactly embedded in $L^q(D; L^r(\omega))$.
\item [{\rm (ii)}] If $q \ge \max\{q_0, q_1\}$ is such that  
$$
\left(1 - \frac{1}{r_0}\right)\left(\frac{1}{N} - \frac{1}{q_1} + \frac{1}{q}\right)
> \frac{1}{r_1}\left(\frac{1}{q_0} - \frac{1}{q}\right),
$$
then $W^{q_0, r_0; q_1, r_1}(D, \omega)$ is compactly embedded in $L^q(D; C(\bar\omega))$.
\item [{\rm (iii)}] If
$$
\left(1 - \frac{1}{r_0}\right)\left(\frac{1}{N} - \frac{1}{q_1}\right)
> \frac{1}{r_1 q_0},
$$
then $W^{q_0, r_0; q_1, r_1}(D, \omega)$ is compactly embedded in $C(\bar D\times \bar\omega))$.
\end{itemize}
\end{proposition}

We have bounds independent of $m$ in $L^3(\Omega; L^3(2\pi, 4\pi))$ for $p\om_t$ and
in $L^2(\Omega; L^r(2\pi, 4\pi))$ for $\nabla p\om$ for all $r>1$ by virtue of \eqref{es1}
and \eqref{es2}. Hence, by Proposition \ref{emb}\,(ii), $\{p\om\}$ is a compact sequence
in $L^q(\Omega; C[2\pi, 4\pi])$ for each $q \in [1,6)$. Furthermore, by \eqref{es4} $\{G_R[p\om]_t\}$
is a bounded sequence in $L^3((2\pi, 4\pi); L^2(\Omega))$. By \eqref{lipc} and \eqref{qua}, we can select
a subsequence (still indexed by $m$) in such a way that
\be{con1}
\left\{
\begin{array}{lll}
p\om \to p \ &\mbox{strongly in }\ &L^4(\Omega; C[2\pi, 4\pi]),\\
G_R[p\om] \to G_R[p]  \ &\mbox{strongly in }\ &L^2(\Omega; C[2\pi, 4\pi]),\\
G_R[p\om]_t \to G_R[p]_t  \ &\mbox{weakly in }\ &L^3((2\pi, 4\pi); L^2(\Omega)).
\end{array}
\right.
\ee
All the other terms in \eqref{e1g}--\eqref{e2g} are linear, so that we can pass to the weak limit
and conclude that the system
\ber \nonumber
&&\ippi\io ((u_{tt} + u_t)\Phi(x,t) + \AA\nas u :\nas\Phi(x,t) + p \dive\Phi(x,t))\dd x \dd t\\ \label{e1l}
&& \qquad = \ippi\io f\Phi(x,t)\dd x\dd t\,,\\ \nonumber
&&\ippi\io ((G_R[p]_t - \dive u_t)\Psi(x,t) + \nabla p \nabla\Psi(x,t))\dd x \dd t\\ \label{e2l}
&&\qquad  =\ippi\io h\Psi(x,t)\dd x \dd t + \ippi\ipo \gamma(x)(p^* - p)\Psi(x,t)\dd s(x) \dd t\,,
\eer
is satisfied for all arbitrarily chosen $2\pi$-periodic test functions
$\Phi\in L^2((2\pi, 4\pi); W^{1,2}_0(\Omega;\real^3))$,
$\Psi\in L^2((2\pi, 4\pi); W^{1,2}(\Omega))$, with the regularity
$u,u_t, u_{tt}, \nas u_t, \nabla p, \nabla p_t \in L^2\pp(\Omega)$, $p_t \in L^3\pp(\Omega)$,
$G_R[p]_t \in L^3((2\pi, 4\pi); L^2(\Omega))$. Choosing $\Psi$ with compact support in $\Omega$,
we obtain from \eqref{es3}, \eqref{es4} that
$$
\|\Delta p\|_{2,\Omega,2\pi} \le C\delta\,,
$$
and similarly
$$
\|\dive \AA \nas u\|_{2,\Omega,2\pi} \le C\delta\,,
$$
so that the identities
\ber  \label{e1s}
u_{tt} + u_t - \dive \AA\nas u &=& \nabla p + f\,,\\ \label{e2s}
G_R[p]_t - \Delta p &=& \dive u_t + h
\eer
with boundary conditions \eqref{bcb} hold in the sense of $L^2\pp(\Omega)$.

For $s>0$, and a function $v\in L^2\pp(\Omega)$ such that $v_t\in L^2\pp(\Omega)$, we define
$$
v\os(x,t) = \frac{1}{\tau}(v(x,t+\tau) - v(x,t))\,.
$$
We have indeed
\be{lims}
\ippi\io |v\os(x,t)|^2 \dd x\dd t\le \ippi\io |v_t(x,t)|^2 \dd x\dd t\,,
\quad \lim_{s\to 0} \ippi\io |v\os(x,t) - v_t(x,t)|^2 \dd x\dd t = 0\,.
\ee
Put $\hat f = \nabla p + f$. Then $\hat f, \hat f_t \in L^2\pp(\Omega)$
by virtue of \eqref{es2} and the hypotheses on $f$. From \eqref{e1s} it follows for $s>0$ that
\be{e6}
u\os_{tt} + u\os_t - \dive \AA\nas u\os = \hat f\os\,.
\ee
We test \eqref{e6} by $u\os_t$ and obtain for a.~e. $t \in \real$ that
$$
\frac{\dd}{\dd t} \io (|u\os_t|^2 + \AA\nas u\os: \nas u\os)(x,t)\dd x
+ \io |u\os_t|^2(x,t)\dd x \le
\io |\hat f\os|^2(x,t)\dd x\,,
$$
so that
\ber \nonumber
&&\hspace{-14mm}\frac{\dd}{\dd t} \io (|u\os_t|^2 + \AA\nas u\os: \nas u\os)(x,t)\dd x
+ \io (|u\os_t|^2 + \AA\nas u\os: \nas u\os)(x,t)\dd x \\ \label{e7}
&\le& \io (|\hat f\os|^2 + \AA\nas u\os: \nas u\os)(x,t)\dd x\,.
\eer
Put
$$
y(t) = \io (|u\os_t|^2 + \AA\nas u\os: \nas u\os)(x,t)\dd x\,, \quad
\beta(t) = \io (|\hat f\os|^2 + \AA\nas u\os: \nas u\os)(x,t)\dd x\,.
$$
Both $y$ and $\beta$ belong to $L^1\pp$ by \eqref{es2}--\eqref{es3} and
\be{e8}
\ippi (y(t) + \beta(t))\dd t \le C\delta^2\,. 
\ee
Then \eqref{e7} is an inequality of the form
$$
\frac{\dd}{\dd t} y(t) + y(t) \le \beta (t)\,,
$$
which implies by the Gronwall argument that
$$
y(t) \le C\ippi (y(\tau)+ \beta(\tau))\dd\tau \le C\delta^2 \quad \mbox{for a.~e. } t\in \real\,,
$$
and we conclude by passing to the limit as $s\to 0$ and using \eqref{lims}
that $u_{tt}, \nas u_t$ belong to $L^\infty(\real; L^2(\Omega))$ with
\be{e9}
\supess_{t\in \real}\io (|u_{tt}|^2 + |\nas u_t|^2)(x,t)\dd x \le C\delta^2\,.
\ee
We see that in \eqref{e2s}, we have $h \in L^\infty(\real; L^2(\Omega))$ by hypothesis,
$\dive u_t \in L^\infty(\real; L^2(\Omega))$ by Korn's inequality and \eqref{e9}, and
$G_R[p]_t \in L^3((2\pi, 4\pi); L^2(\Omega))$ by \eqref{con1}, with bounds proportional
to $\delta$. Hence, with the notation \eqref{lrq},
\be{e10}
\|\Delta p\|_{2,3} \le C\delta\,.
\ee
By hypotheses about the regularity of $\partial\Omega$ and of the boundary data, and by \eqref{es2},
we have
\be{e11}
\left\|\frac{\partial^2 p}{\partial x_i \partial x_j}\right\|_{2,3} \le C\delta\,,\quad
\left\|\frac{\partial^2 p}{\partial x_i \partial t}\right\|_{2,2} \le C\delta
\ee
for all $i,j = 1,2,3$. From Proposition \ref{emb}\,(i) we obtain that
\be{e12}
\|\nabla p\|_{r,r} \le C\delta
\ee
for $r < 14/3$. Combining this result with the fact that
\be{e13}
\|p_t\|_{3,3} \le C\delta
\ee
which follows from \eqref{es2}, we can use Proposition \ref{emb}\,(iii) and conclude that
\be{e14}
\max\{|p(x,t)|: (x,t) \in \bar\Omega \times [2\pi. 4\pi]\} \le C\delta
\ee
provided $r > 9/2$. Thus, choosing
$$
r \in \left(\frac92, \frac{14}3\right),
$$ 
we see that if $\delta>0$ is chosen sufficiently small, then $|p(x,t)|$ does not exceed the critical
value $R$, and from \eqref{conve} we infer that the solution of \eqref{e1s}--\eqref{e2s} that we have
constructed is the desired solution of \eqref{e1b}--\eqref{e2b}, which we wanted to prove.


\section*{Conclusion}

We have proved that a model for the propagation of periodic mechanical waves
inside an elastic partially saturated porous body with capillary hysteresis represented
by a Preisach operator is well posed under periodic boundary forcing provided
the boundary forces are sufficiently small. The meaning of the smallness condition 
is to keep the pressure values within the convexity domain of the Preisach operator
and exploit the hysteresis second order energy inequality.


{\small
\begin{thebibliography}{99}

\bibitem{a2006} Albers B. Monochromatic surface waves at the interface between poroelastic and fluid
halfspaces. {\em Proc. Royal Soc. A}. 2006; {\bf 462}: 701--723.  


\bibitem{AM} Albers B. Modeling the hysteretic behavior of the capillary pressure in partially saturated porous media -- a review. {\em Acta Mechanica}. 2014; {\bf 225} (8): 2163--2189. 
 
\bibitem{CEER} Albers B. Main drying and wetting curves of soils -- on measurements, prediction and influence on wave propagation. {\em Engineering Transactions}. 2015; {\bf 63} (1): 5--34.

\bibitem{Murphys} Albers B. On the influence of the hysteretic behavior of the capillary pressure on the wave propagation in partially saturated soils. To appear in {\em Journal of Physics: Conference Series, Proceedings of the 7th International Workshop on Multi-Rate Processes and Hysteresis}. 2014.

\bibitem{Hab} Albers B. \emph{Modeling and Numerical Analysis of Wave
Propagation in Saturated and Partially Saturated Porous Media}. Habilitation
Thesis. Ver\"{o}ffentlichungen des Grundbauinstitutes der Technischen
Universit\"{a}t Berlin, Shaker: Aachen; 2010.

\bibitem{ak} Albers B., Krej\v{c}\'{\i} P. Unsaturated porous media flow with
thermomechanical interaction. Accepted for publication in MMAS.

\bibitem{ak2016} Albers B., Krej\v{c}\'{\i} P.
Hysteresis in unsaturated porous media -- two models for wave propagation and engineering applications. In
{\em Continuous Media with Microstructure 2}, Eds. B.~Albers and M.~Kuczma, Springer, 2016: 217--229.

\bibitem{bv1} Bagagiolo F, Visintin A. Hysteresis in filtration through porous media.
{\em Z. Anal. Anwendungen}. 2000; {\bf 19} (4): 977--997.

\bibitem{bv2} Bagagiolo F, Visintin A. Porous media filtration with hysteresis.
{\em Adv. Math. Sci. Appl.} 2004; {\bf 14} (2): 379--403.

\bibitem{bin} Besov OV, Il'in VP, Nikol'ski\u{\i} SM.
  \emph{Integral Representations of Functions and Imbedding Theorems}.
  Scripta Series in Mathematics.
  Halsted Press (John Wiley \& Sons):
  New York-Toronto, Ont.-London; 1978 (Vol. I), 1979 (Vol. II). Russian
version Nauka: Moscow; 1975.

\bibitem{CosentiniFoti} Cosentini RM, Foti S. Evaluation of porosity and degree of saturation from
seismic and electrical data. {\em G\'eotechnique} 2014; {\bf 64} (4): 278--286.

\bibitem {mhd} Eleuteri M, Kopfov\'a J, Krej\v{c}\'{\i} P.
Magnetohydrodynamic flow with hysteresis. {\em SIAM J. Math. Anal.} 2009; {\bf 41}: 435--464. 

\bibitem{Flynn} Flynn D. {\em Modelling the flow of water through multiphase porous media with the Preisach model}. PhD Thesis. University College Cork; 2008.

\bibitem{flpokr} Flynn D, McNamara H, O'Kane JP, Pokrovski\u{\i} AV.
Application of the Preisach model to soil-moisture hysteresis.
{\em The Science of Hysteresis, Volume 3}, Eds. Bertotti G, Mayergoyz I. Academic Press: Oxford; 2006: 689--744.

\bibitem{Foti} Foti S, Lai CG, Rix, GJ, Strobbia C. {\em Surface Wave Methods for Near-Surface Site Characterization.}
Boca Raton: CRC Press; 2014.

\bibitem{Gan} Gan WS. {\em Acoustical Imaging: Techniques and Applications for Engineers.} Wiley; 2012.

\bibitem{hrrp} Haverkamp R, Reggiani P, Ross PJ, Parlange J-Y.
Soil water hysteresis prediction model based on theory and geometric scaling. \emph{Environmental Mechanics,
Water, Mass and Engergy Transfer in the Biosphere}. Eds. Raats PAC, Smiles D, Warrick AW.
American Geophysical Union: 2002; 213--246.

\bibitem{kp} Krasnosel'ski\u{\i} MA, Pokrovski\u{\i} AV.
{\em Systems with Hysteresis}. Springer: Berlin; 1989. 
Russian edition: Nauka: Moscow; 1983.

\bibitem{Kraut} Krautkr\"amer J, Krautkr\"amer H. {\em Ultrasonic Testing of Materials.}
Springer-Verlag, Berlin; 1990.

\bibitem{book} Krej\v c\'{\i} P. {\em Hysteresis, Convexity and Dissipation in
Hyperbolic Equations\/}. Gakuto Intern. Ser. Math. Sci. Appl. Vol.~{\bf 8}.
Gakkot\={o}sho: Tokyo; 1996.

\bibitem {max} Krej\v{c}\'{\i} P. On Maxwell equations
with the Preisach hysteresis operator: the one-dimensional
time-periodic case. {\em Apl. Mat.} 1989; {\bf 34}: 364--374.

\bibitem {pr} Preisach F. \"Uber die magnetische Nachwirkung.
{\em Z.~Phys.} 1935; {\bf 94}: 277--302 (in German).

\bibitem {Roesset} Roesset JM. Nondestructive dynamic testing of soils and pavements.
{\em Tamkang Journal of Science and Engineering} 1998; 61--81.

\bibitem{Shull} Shull PJ. {\em Nondestructive Evaluation: Theory, Techniques, and Applications}.
Marcel Dekker Inc.; 2002.

\bibitem{Veenstra} Veenstra M. White DJ, Schaefer VR.
{\em Synthesis of Nondestructive Testing Technologies for Geomaterial Applications}. Iowa State University; 2005.



\end{thebibliography}
}
\end{document}